\documentclass[12pt]{amsart}
\usepackage{amsfonts}
\usepackage{amsmath}
\usepackage{amssymb}
\usepackage{a4}

\newcommand{\heute}{13 August 2008}

\theoremstyle{plain}
\newtheorem{theorem}{Theorem}[section]

\newtheorem{lemma}[theorem]{Lemma}

\newtheorem{proposition}[theorem]{Proposition}
\theoremstyle{remark}
\newtheorem{conj}[theorem]{Conjecture}
\newtheorem{remark}[theorem]{Remark}
\newtheorem{ques}[theorem]{Question}
\newtheorem*{defn}{Definition}
\newtheorem*{rk}{Remark}

\newcommand{\dashTwo}[1]{\textup{(\ref{two}${}'$)}}

\newcommand{\ignore}[1]{}

\newcommand{\f}[1][p]{\mathbb{F}_{#1}}

\newcommand{\Gro}[1]{Gr\"ob\-ner}

\newcommand{\depth}{\operatorname{depth}}

\newcommand{\Ann}{\operatorname{Ann}}

\newcommand{\Reg}{\operatorname{Reg}}

\newcommand{\CMd}{\delta}
\newcommand{\gtD}{\operatorname{gtD}}
\newcommand{\De}{e}

\DeclareMathSymbol\normal{\mathrel}{AMSa}{"43}

\newcommand{\prank}[1][p]{\text{$#1$-rk}}

\begin{document}

\title{Testing {B}enson's regularity conjecture}
\author{David J. Green}
\address{Mathematical Institute \\
University of Jena \\ D-07737 Jena \\
Germany}
\email{green@minet.uni-jena.de}
\thanks{Travel assistance from DFG grant GR 1585/4-1}
\date{\heute}

\begin{abstract}
\noindent
D.~J. Benson conjectures that the Castelnuovo-Mumford regularity of the
cohomology ring of a finite group is always zero.
More generally he conjectures that there is always a
\emph{very strongly quasi-regular} system of parameters.
Computer calculations show that the second conjecture holds
for all groups of order less than $256$.
\end{abstract}

\maketitle
\newlength{\djglength}

\section{Introduction}
\noindent
This paper is concerned with a conjecture of
D.~J. Benson~\cite{Benson:DicksonCompCoho} about the commutative algebra
of group cohomology rings.
There are several results relating the group structure
of a finite group~$G$ to the commutative algebra of its cohomology ring
$H^*(G) = H^*(G,k)$ with coefficients in a field of characteristic~$p$.
For the Krull dimension and depth we have the following inequalities,
where $S$ denotes a Sylow $p$-subgroup of~$G$.
Recall that the $p$-rank of $G$ is the rank of the largest elementary
abelian $p$-subgroup.
\begin{equation}
\label{eqn:QuillenDuflot}
\prank(Z(S)) \leq \depth H^*(S) \leq \depth H^*(G)
\leq \dim H^*(G) = \prank(G) \, .
\end{equation}
See Evens' book~\cite{Evens:book} for proofs of the
first inequality (Duflot's theorem) and the last one (due to Quillen)\@.
The second inequality is Theorem~2.1 of Benson's paper~\cite{Benson:NYJM2}
and
``must be well known''\@. The third is
automatic for finitely generated connected $k$-algebras.
Note that the dimension and depth only depend on $G,p$: not on~$k$.
These inequalities motivate the following definitions.

\begin{defn}
Let $G, p, k,S$ be as above. The
group-theoretic defect $\gtD_p(G)$,
and the
Cohen--Macaulay defect $\CMd_p(G)$
are defined by
\begin{xalignat*}{2}
\gtD(G) & = \prank(G) - \prank(Z(S))
&
\CMd(G) & = \dim H^*(G,k) - \depth H^*(G,k)
\end{xalignat*}
\end{defn}

\noindent
It follows from Eqn.~\eqref{eqn:QuillenDuflot} that
\begin{xalignat}{3}
\label{eqn:gtD-CMd}
0 \leq \CMd(G) & \leq \gtD(G) &
\gtD(G) & = \gtD(S) & \CMd(G) & \leq \CMd(S) \, .
\end{xalignat}
The term Cohen--Macaulay defect (sometimes deficiency) is already in
use among workers in the field.
%
To state Benson's conjectures we need some terminology.

\begin{defn}
Let $p,k$ be as above. Let $A$ be a graded commutative $k$-algebra which
is both connected and finitely generated. Connected means that $A^0=k$
and $A^{{}<0}=0$.
Let $\zeta_1,\ldots,\zeta_r$ be a system of homogeneous elements in~$A^{{}>0}$,
and set $n_i = \left|\zeta_i\right| > 0$.
\begin{enumerate}
\item
The system is called a filter-regular system of parameters if multiplication
by~$\zeta_{i+1}$ has finite-dimensional kernel as an endomorphism of
$A/(\zeta_1,\ldots,\zeta_i)$ for each $0 \leq i \leq r$,
where $\zeta_{r+1}=0$.
Observe that a filter-regular system of parameters really is a system of
parameters.
\item
A very strongly quasi-regular system of parameters is a system which is a
filter-regular system of parameters by virtue of the property that
this kernel is restricted to degrees${} \leq n_1+\cdots+n_i+d_i$
for each $0 \leq i \leq r$, where $d_r=-r$ and $d_i = -i-1$ for all $i < r$.
\end{enumerate}
\end{defn}

\begin{theorem}[Benson]
\label{thm:portfolio}
Let $G$ be a finite group, $p$ a prime number and $k$ a field
of characteristic~$p$.
\begin{enumerate}
\item
\label{enum:DicksonFilter}
The cohomology ring $A = H^*(G,k)$ does have filter-regular systems of
parameters: the Dickson invariants (suitably interpreted) form one.
\item
\label{enum:FilterEqual}
Either every filter-regular system of parameters in~$A$ is very strongly
quasi-regular, or none are.
\item
\label{enum:OkuSas}
If the Cohen--Macaulay defect of~$G$ satisfies $\CMd(G) \leq 2$ then
there is a very stong quasi-regular system of parameters in~$A$.
In particular $\CMd(G) \leq 2$ holds for all 267 groups of order~$64$.
\item
\label{enum:hierarchy}
If there is a very strongly quasi-regular system of parameters in~$A$,
then the Castelnuovo--Mumford regularity of~$A$ is zero.
\end{enumerate}
\end{theorem}

\begin{proof}
The main reference is Benson's paper~\cite{Benson:DicksonCompCoho}\@.
Part~\ref{enum:DicksonFilter}) is Coroll.~9.8\@ and
Part~\ref{enum:FilterEqual}) is Coroll.~4.7(c)\@, whereas
Part~\ref{enum:hierarchy}) follows from Coroll.~4.7(c) and Theorem~4.2\@.
The first statement of Part~\ref{enum:OkuSas}) is Theorem~1.5
of~\cite{Benson:DicksonCompCoho}\@; 
the second one was observed by Carlson~\cite{Carlson:Online3,CarlsonTownsley},
who computed the cohomology ring of every group of order~$64$\@.
The reader may find the tabulated data in~\cite[Appendix]{Benson:MSRI} useful.
\end{proof}

\begin{rk}
A weaker version of the $\CMd(G)=2$ case of~\ref{enum:OkuSas}) was also
proved by Okuyama and Sasaki. It is a shame that their paper~\cite{OkuSas:Hsop}
appeared so late: I know that it had completed the refereeing process by the
end of April 2001, but it had been superseded by the time it was finally
published in 2004\@.
\end{rk}

\begin{conj}[Benson~\cite{Benson:DicksonCompCoho}]
\label{conj:Reg}
Let $G$ be a finite group, $p$ a prime number and $k$ a field
of characteristic~$p$. The cohomology ring $H^*(G,k)$ has Castelnuovo--Mumford
regularity zero.
\end{conj}

\begin{conj}
\label{conj:VSQR}
Let $G$ be a finite group, $p$ a prime number and $k$ a field
of characteristic~$p$. The cohomology ring $H^*(G,k)$ always contains a
very strongly quasi-regular system of parameters.
\end{conj}

\begin{rk}
Conjecture~\ref{conj:Reg} is Benson's Conjecture~1.1\@.
By Theorem~\ref{thm:portfolio}~\ref{enum:hierarchy}) it is a weak form of
Conjecture~\ref{conj:VSQR}, which is only implicitly present in Benson's paper.
Kuhn has shown that Conjecture~\ref{conj:Reg} has applications to
the study of central essential cohomology~\cite{Kuhn:Cess}\@.
\end{rk}

\noindent
The conjectures have been verified in two families of cases.
Benson showed in~\cite{Benson:wreathReg} that if Conjecture~\ref{conj:Reg}
holds for~$H$ then it also holds for the wreath product $G = H \wr C_p$.
And the second verification is the following theorem, the main
result of the present paper.

\begin{theorem}
\label{thm:main}
Conjecture~\ref{conj:VSQR} holds for every group of order less than 256\@.
\end{theorem}

\begin{proof}
By Theorem~\ref{thm:portfolio}~\ref{enum:OkuSas}) a counterexample has to
have $\CMd(G) \geq 3$. By Proposition~\ref{prop:only128} the only groups
of order less than $256$ satisfying $\CMd(G) \geq 3$ have order~$128$
and satisfy $\CMd(G)=3$. By Proposition~\ref{prop:calc} there are fourteen
groups of order~$128$ with $\CMd(G)=3$, and each of these satisfies
the conjecture.
\end{proof}

\section{Reduction to the case $|G|=128$}

\begin{proposition}
\label{prop:only128}
Let $G$ be a group of order less than $256$. Then $\CMd(G) \leq 3$;
and if $\CMd(G)=3$ then $|G|=128$.
\end{proposition}

\begin{proof}
Let $S$ be a Sylow $p$-subgroup of~$G$. In view of the inequality
$\CMd(G) \leq \CMd(S)$ and the restriction $|G| < 256$
it suffices to consider the case where $G$ itself is a $p$-group.

So suppose $G$ is a $p$-group with $\CMd(G) \geq 3$.
By Lemma~\ref{lemma:Jordan} below it follows that $p=2$,
that $\CMd(G)=3$, and that either $|G|=64$ or $|G|=128$.
But Carlson's computations [see Theorem~\ref{thm:portfolio}\ref{enum:OkuSas})
above] show that $\CMd(G) \leq 2$ if $\left|G\right|=64$.
\end{proof}

\begin{lemma}
\label{lemma:Jordan}
Let $G$ be a finite group and $p$ a prime number.
\begin{enumerate}
\item
\label{enumi:Jordan3}
If $\CMd(G) \geq 3$
or more generally if $\gtD(G) \geq 3$ then $p^5$ divides the order of~$G$.
If $p=2$ or $p=3$ then $p^6$ must divide $\left|G\right|$.
\item
\label{enumi:Jordan4}
If $\CMd(G) \geq 4$ or more generally if $\gtD(G) \geq 4$ then $p^6$ divides
the order of~$G$.  If $p=2$ or $p=3$ then $p^7$ must divide $\left|G\right|$.
\item
\label{enumi:Jordan128}
If $p=2$ and $\gtD(G) \geq 4$ then $\left|G\right|$ is divisible by $256$.
\end{enumerate}
\end{lemma}

\begin{proof}
\ref{enumi:Jordan3}):
By Eqn.~\eqref{eqn:gtD-CMd} we have $\CMd(G) \leq \gtD(G)$.
It is apparent from the definition that a finite group and its Sylow
$p$-subgroups have the same group-theoretic defect.
So it suffices to consider the case where $G$~is a $p$-group and
$\gtD(G)\geq 3$.

Every nontrivial $p$-group has a centre of rank at least one. So a $p$-group
with $\gtD \geq 3$ must have a subgroup which is elementary abelian of rank~$4$.
It must be nonabelian too, so the order must be at least~$p^5$.

Suppose that
there is such a group of order~$p^5$. Then the centre is cyclic of order~$p$,
and there is an elementary abelian subgroup~$V$ of order~$p^4$.
This $V$~is a maximal subgroup of a $p$-group and therefore normal.
Let $a \in G$ lie outside~$V$. Then $G = \langle a,V\rangle$
and the conjugation action of $a$~on $V$ must be nontrivial of order~$p$.
So the minimal polynomial of the action divides $X^p-1 = (X-1)^p$.
This means that the action has a Jordan normal form with sole eigenvalue~$1$.
The eigenvectors in each Jordan block belong to the centre of~$G$. So as
the centre is cyclic there can only be one Jordan block, of size~$4$\@.
But for $p=2$
the size~$3$ Jordan block does not square to the identity, so there can
be no blocks of size $3$~or higher. Similarly there can be no size~$4$ block
for $p=3$, since it does not cube to the identity. So for $p=2,3$ there must
be more than one Jordan block.
\medskip

\noindent
\ref{enumi:Jordan4}): analogous.
\medskip

\noindent
\ref{enumi:Jordan128}): We need to show that if $\left|G\right|=128$
then $\gtD(G) \geq 4$ cannot happen. Such a group would need to contain
an elementary abelian $2$-group~$V$ satisfying the equation
$\prank(V) = 4 + \prank(Z(G))$.
The rank of the centre must be at least one, and if it is three or more
then $V=G$ and therefore $G$ is abelian, a contradiction. If the centre has
rank two, then $V$~has index two in~$G$ and we are in the same situation
as in~\ref{enumi:Jordan3}), except now we want two Jordan blocks. But we
need at least three, since $V$ has rank six and only Jordan blocks
of size one or two are allowed.

If the centre has rank one, then $V$~has rank five and we may pick a
subgroup~$H$ with $V \leq H \leq G$ and $[G:H]=[H:V]=2$. If $H$ is elementary
abelian then we are back in the immediately preceding case of an order
two action on a rank six elementary abelian. If~$H$ is not elementary abelian
then the usual Jordan block considerations mean that $C$ has rank at least~$3$,
where $C$~is the unique largest central elementary abelian of~$H$. By
applying the Jordan block considerations to the conjugation action
of $G/H$ on~$C$, we see that $Z(G)$ must have rank at least two, a
contradiction.
\end{proof}

\begin{remark}
In fact there are only two groups of
order $64$ with $\gtD=3$.
Here are their numbers in the
Hall--Senior list~\cite{HallSenior} and in the Small Groups
Library~\cite{BeEiOBr:Millennium}\@.
Their defects are taken from the tables in~\cite[Appendix]{Benson:MSRI}\@.
\settowidth{\djglength}{$000$}
\[ \begin{array}{c|c|c}
\text{Small Group} & \text{Hall--Senior} & \CMd(G) \\
\hline
\makebox[\djglength][r]{$32$} & \makebox[\djglength][r]{$250$} & 2 \\
\makebox[\djglength][r]{$138$} & \makebox[\djglength][r]{$259$} & 1
\end{array} \]
As we shall see below there are  $14$ groups of order~$128$
with $\CMd(G)=3$.
\end{remark}

\begin{remark}
For $p \geq 5$ let $G$ be the following semidirect product group of order~$p^5$:
there is a rank four elementary abelian on the bottom and a cyclic group
of order~$p$ on top. The conjugation action is a size 4 Jordan block. 
This group has $\gtD(G)=3$, since
\[
\begin{pmatrix}
1 & 1 & 0 & 0 \\ 0 & 1 & 1 & 0 \\ 0 & 0 & 1 & 1 \\ 0 & 0 & 0 & 1 \end{pmatrix}^n
= \begin{pmatrix} 1 & \binom n1 & \binom n2 & \binom n3 \\
0 & 1 & \binom n1 & \binom n2 \\ 0 & 0 & 1 & \binom n1 \\ 0 & 0 & 0 & 1
\end{pmatrix} \, .
\]
\end{remark}

\section{The groups of order 128}
\noindent
Let $G$ be a group of order~$128$.
By Lemma~\ref{lemma:Jordan}\ref{enumi:Jordan128}) we have $\gtD(G) \leq 3$.

\begin{proposition}
\label{prop:calc}
Only $57$ out of the $2328$ groups of order $128$ satisfy $\gtD(G) = 3$.
Of these $57$ groups, $43$ satisfy $\CMd(G) \leq 2$. The remaining $14$ groups
satisfy $\CMd(G)=3$. Each of these $14$ groups of order $128$ with
$\CMd(G)=3$ satisfies Conjecture~\ref{conj:VSQR}\@.

According to the numbering of the Small Groups
Library~\cite{BeEiOBr:Millennium} these fourteen groups are:
numbers
$36$, $48$, $52$, $194$, $515$, $551$, $560$, $561$, $761$, $780$, $801$,
$813$, $823$ and $836$.
\end{proposition}

\begin{proof}
By machine computation.
Inspecting the Small Groups library,
one sees that there are 57 groups with $\gtD(G)=3$. See
Appendix~\ref{sect:appendix} for a discussion of how the $p$-rank
is computed.

\begin{table}
\settowidth{\djglength}{$0000$}
\newcommand{\gr}[1]{\makebox[\djglength][r]{$#1$}}
\newcommand{\ugr}[1]{\makebox[\djglength][r]{\underline{$#1$}}}
$\begin{array}{|c|cccc|c|c|cccc|}
\cline{1-5} \cline{7-11}
\text{gp} & K & d & r & \CMd & \quad & \text{gp} & K & d & r & \CMd \\
\cline{1-5} \cline{7-11}
\ugr{36} & 5 & 2 & 2 & 3 & & \gr{850} & 5 & 3 & 2 & 2 \\
\ugr{48} & 5 & 2 & 2 & 3 & & \gr{852} & 4 & 2 & 1 & 2 \\
\ugr{52} & 4 & 1 & 1 & 3 & & \gr{853} & 4 & 2 & 1 & 2 \\
\ugr{194} & 5 & 2 & 2 & 3 & & \gr{854} & 4 & 2 & 1 & 2 \\
\gr{513} & 5 & 3 & 2 & 2 & & \gr{859} & 4 & 2 & 1 & 2 \\
\ugr{515} & 5 & 2 & 2 & 3 & & \gr{860} & 4 & 2 & 1 & 2 \\
\gr{527} & 4 & 2 & 1 & 2 & & \gr{866} & 4 & 2 & 1 & 2 \\
\ugr{551} & 5 & 2 & 2 & 3 & & \gr{928} & 4 & 3 & 1 & 1 \\
\ugr{560} & 4 & 1 & 1 & 3 & & \gr{929} & 4 & 2 & 1 & 2 \\
\ugr{561} & 4 & 1 & 1 & 3 & & \gr{931} & 4 & 2 & 1 & 2 \\
\gr{621} & 5 & 3 & 2 & 2 & & \gr{932} & 4 & 2 & 1 & 2 \\
\gr{623} & 4 & 2 & 1 & 2 & & \gr{934} & 4 & 2 & 1 & 2 \\
\gr{630} & 5 & 3 & 2 & 2 & & \gr{1578} & 6 & 4 & 3 & 2 \\
\gr{635} & 4 & 2 & 1 & 2 & & \gr{1615} & 4 & 3 & 1 & 1 \\
\gr{636} & 4 & 2 & 1 & 2 & & \gr{1620} & 4 & 2 & 1 & 2 \\
\gr{642} & 4 & 2 & 1 & 2 & & \gr{1735} & 5 & 3 & 2 & 2 \\
\gr{643} & 4 & 3 & 1 & 1 & & \gr{1751} & 4 & 2 & 1 & 2 \\
\gr{645} & 4 & 3 & 1 & 1 & & \gr{1753} & 4 & 3 & 1 & 1 \\
\gr{646} & 4 & 2 & 1 & 2 & & \gr{1755} & 5 & 4 & 2 & 1 \\
\gr{740} & 4 & 2 & 1 & 2 & & \gr{1757} & 4 & 3 & 1 & 1 \\
\gr{742} & 4 & 2 & 1 & 2 & & \gr{1758} & 4 & 3 & 1 & 1 \\
\gr{753} & 5 & 3 & 2 & 2 & & \gr{1759} & 4 & 3 & 1 & 1 \\
\ugr{761} & 5 & 2 & 2 & 3 & & \gr{1800} & 4 & 2 & 1 & 2 \\
\gr{764} & 4 & 2 & 1 & 2 & & \gr{2216} & 5 & 4 & 2 & 1 \\
\ugr{780} & 4 & 1 & 1 & 3 & & \gr{2222} & 5 & 3 & 2 & 2 \\
\ugr{801} & 4 & 1 & 1 & 3 & & \gr{2264} & 5 & 3 & 2 & 2 \\
\ugr{813} & 4 & 1 & 1 & 3 & & \gr{2317} & 4 & 3 & 1 & 1 \\
\ugr{823} & 4 & 1 & 1 & 3 & & \gr{2326} & 4 & 4 & 1 & 0 \\
\cline{7-11}
\ugr{836} & 4 & 1 & 1 & 3 \\
\cline{1-5}
\end{array}$
\vspace*{5pt}
\caption{For each group of order $128$ with $\gtD(G)=3$, we give its number
in the Small Groups library, the Krull dimension $K$ and
depth $d$ of $H^*(G)$, the rank~$r = K - 3$ of $Z(G)$ and the
Cohen--Macaulay defect $\CMd = K - d$. Underlined entries have $\CMd = 3$\@.
Notation based on that of~\cite[Appendix]{Benson:MSRI}\@.}
\label{table:128}
\end{table}

These $57$ groups are listed
in Table~\ref{table:128}\@. The cohomology rings of these $57$ groups were
computed using an improved version of the author's cohomology
program~\cite{habil}\@.
These cohomology rings may be viewed online~\cite{fullBens}\@.
The cohomology rings were calculated using Benson's test for
completion~\cite[Thm 10.1]{Benson:DicksonCompCoho}\@.

Benson's test involves constructing a filter-regular system of parameters
and determining in which degrees it is not strictly regular. This means that
one automatically determines whether the group satisfies
Conjecture~\ref{conj:VSQR} when one computes cohomology using Benson's test.
The Cohen--Macaulay defect is another by-product of a computation
based on Benson's test.

The value of $\CMd(G)$ for each of the $57$ groups
is given in Table~\ref{table:128}\@. The fourteen groups listed in the
statement of the proposition are indeed
the only ones with $\CMd(G)=3$.
The computations showed that these $14$ groups do satisfy the conjecture.
\end{proof}

\begin{rk}
The test is phrased in such a way that it is easy to implement. With one
exception: it is not immediately clear how to construct a filter-regular
system of parameters in low degrees. This point is discussed in the next
section.
\end{rk}

\begin{rk}
The computation that took the longest time was group number $836$,
one of the $\CMd(G)=3$ groups. Its cohomology ring has $65$ generators
and $1859$ generators.
\end{rk}

\begin{rk}
The distribution of these 57 groups by Cohen--Macaulay defect is as
follows:
\[
\begin{array}{c|cccc}
\CMd(G) & 0 & 1 & 2 & 3 \\
\hline
\#G & 1 & 11 & 31 & 14
\end{array}
\]
\end{rk}

\begin{rk}
Some of these groups have been studied before.
Groups 928 and 1578 are wreath products: $D_8 \wr 2$ and $2^3 \wr 2$
respectively. By the Carlson--Henn result~\cite{CaHe:Wreath} one has
$\CMd(D_8 \wr 2)=1$ and $\CMd(2^3\wr2)=2$.
Groups 850 of order 128 is a direct product of the form $G=H \times 2$,
where $H$ is group number 32 of order $64$; and the same applies to group
1755 of order 128 and group 138 of order $64$\@. It is immediate that
$\CMd(G)=\CMd(H)$ for such groups, so Carlson's work guarantees $\CMd(G)\leq2$
for both these groups of order~$128$.

Group number $2326$ is the extraspecial group~$2^{1+6}_+$; Quillen showed that
its cohomology ring is Cohen--Macaulay~\cite{Quillen:Extraspecial}\@.
Group number $931$ is the Sylow $2$-subgroup of the Mathieu groups
$M_{22}$~and $M_{23}$; its cohomology was studied by Adem and
Milgram~\cite{AdMi:M22}\@.
Group number $934$ is the Sylow $2$-subgroup of the Janko group~$J_2$; its
cohomology ring was calculated by Maginnis~\cite{Maginnis:J2}.

I am not aware of any previous cohomological investigations concerning
the other two groups that I can name. One of these is group $932$, the
Sylow $2$-subgroup of $G_2(3):2$. The other is group
number $836$, the Sylow $2$-subgroup of one double cover of the Suzuki
group $\mathit{Sz}(8)$. This group (number 836) turned out to have the most
complicated cohomology ring in the study.
\end{rk}

\section{The weak rank-restriction condition}
\noindent
How does one construct a filter-regular system of parameters in a cohomology
ring which is defined over the prime field~$\f$ and also lies in low degree?
An efficient implementation of Benson's test calls for an answer to this
question.

Benson shows that the Dickson invariants (suitably interpreted) form a
filter-regular system of parameters. This means a sequence of cohomology
classes in $H^*(G)$ whose restrictions to the elementary abelian subgroups
of~$G$ are (powers of) the appropriate Dickson invariants. Given information
about restriction to subrings it is a straightforward task to compute classes
with the appropriate restriction patterns. However the degrees involved
can be large.

\begin{defn}
(c.f.\@ \cite[\S8]{Benson:DicksonCompCoho}) \quad
Let $G$ be a $p$-group with $\prank(G)=K$. Let
$C = \Omega_1(Z(G))$. Homogeneous elements $\zeta_1,\ldots,\zeta_K
\in H^*(G)$ satisfy the \emph{weak rank restriction condition} if
for each rank elementary abelian subgroup $V \geq C$ of $G$ the
following holds, where $s = \prank(V)$:
\begin{quote}
\noindent
The restrictions of $\zeta_1,\ldots,\zeta_s$ to~$V$ form a
homogeneous system of parameters for $H^*(V)$; and
the restrictions of $\zeta_{s+1},\ldots,\zeta_K$ are zero.
\end{quote}
\end{defn}

\begin{lemma}
If $\zeta_1,\ldots,\zeta_K \in H^*(G)$ satisfy the weak rank restriction
condition then they constitute a filter-regular system of parameters.
\end{lemma}

\begin{proof}
By a well known theorem of Quillen (see e.g.\@ Evens' book~\cite{Evens:book}),
$\zeta_1,\ldots,\zeta_K$ is a homogeneous system of parameters for $H^*(G)$.
The proof of Theorem~9.6 of~\cite{Benson:DicksonCompCoho} applies
just as well to parameters satisfying the weak rank restriction condition,
because if $E$ is an arbitrary rank~$i$ elementary abelian subgroup of~$G$,
then setting $V = \langle C, E\rangle$ one has $V \geq C$ and
$C_G(V)=C_G(E)$, yet the rank of $V$ is at least as large as the rank of~$E$.
So the restrictions of $\zeta_1,\ldots,\zeta_i$ to $H^*(C_G(E))$ do form a
regular sequence, by the same argument based on the Broto--Henn approach
to Duflot's theorem.
\end{proof}

\begin{lemma}
\label{lemma:filterConstruct}
Let $G$ be a $p$-group with $\prank(G)=K$ and $\prank(Z(G))=r$. Let
$C = \Omega_1(Z(G))$. Suppose that homogeneous elements
$\zeta_1,\ldots,\zeta_K \in H^*(G)$ satisfy the following conditions:
\begin{enumerate}
\item
The restrictions of $\zeta_1,\ldots,\zeta_r$ to $H^*(C)$
form a regular sequence there; and
\item
For each rank $r+s$ elementary abelian subgroup $V \geq C$ of $G$ the
restrictions of $\zeta_{r+s+1},\ldots,\zeta_K$ to $V$ are zero,
and for $1 \leq i \leq s$ the restrictions of $\zeta_{r+i}$ to~$V$ is
a power of the $i$th Dickson invariant in $H^*(V/C)$.
\end{enumerate}
Then $\zeta_1,\ldots,\zeta_K \in H^*(G)$ is a filter-regular system of
parameters for $H^*(G)$. Moreover such systems of parameters exist.
\end{lemma}

\begin{rk}
By the $i$th Dickson invariant I mean the one which restricts nontrivially
to dimension $i$ subspaces, but has zero restriction to smaller subspaces.
That is, if $i < j$ then the $i$th Dickson invariant lies in lower degree
than the $j$th Dickson invariant.
\end{rk}

\begin{proof}
Such classes clearly satisfy the weak rank restriction condition. The existence
of $\zeta_1,\ldots,\zeta_r$ already follows from Evens' theorem that
the cohomology ring of an arbitrary subgroup $H \leq G$ is a finitely
generated module over the image of restriction from~$G$.
For the $\zeta_{r+i}$: these are given
by restrictions to each elementary abelian subgroup, and these restrictions
satisfy the compatibility conditions that one expects from genuine
restrictions, c.f.\@ Quillen's work on the spectrum of a cohomology
ring~\cite{Quillen:Spectrum}. This means that -- on raising these defining
restrictions by sufficiently high $p$th powers -- the $\zeta_{r+i}$
do indeed exist.
\end{proof}

\begin{rk}
The point is that Lemma~\ref{lemma:filterConstruct} is a recipe for
constructing a filter-regular system of parameters.
Recent work of
Kuhn~\cite{Kuhn:Cess} in fact shows that one can choose the generators of
$H^*(G)$ in such a way that $\zeta_1,\ldots,\zeta_r$ may be chosen from
amongst these generators.

An additional saving follows from the fact that if $\zeta_1,\ldots,\zeta_K$ is
a system of parameters and $\zeta_1,\ldots,\zeta_{K-1}$ is filter-regular,
then the whole system is automatically filter-regular. This means that
one can replace the $\zeta_K$ of Lemma~\ref{lemma:filterConstruct} by any
element that completes a system of parameters.

In earlier calculations, filter-regular parameters were constructed by hand
on a trial and error basis. Subsequently most calculations were performed
or re-performed using the parameter choice method of
Lemma~\ref{lemma:filterConstruct}\@. In the worst case calculations
this meant finishing the computation in degree 17, although the
presentation was finished earlier.
\end{rk}

\section{The $a$-invariants}
Let $k$~be a field and $R$ a connected finitely presented graded commutative
$k$-algebra. Let $M$ be a finitely generated graded $R$-module, and
$\mathfrak{m}$ the ideal in~$R$ of all elements in positive degree.
The $a$-invariants of $M$ are defined by
\[
a_{\mathfrak{m}}^i(M) = \max \{m \mid H^{i,m}_{\mathfrak{m}}(M) \neq 0 \} \, ,
\]
with $a_{\mathfrak{m}}^i(M) = -\infty$ if $H^i_{\mathfrak{m}}(M)=0$.
One can then take
\[
\Reg(M) = \max_{i \geq 0}\{a_{\mathfrak{m}}^i(M) + i \}
\]
as the definition of the Castelnuovo--Mumford regularity of~$M$.
Table~\ref{table:128aInvts} lists the $a$-invariants of $H^*(G)$ for
the $57$ groups of
order 128 with $\gtD(G)=3$\@.

\begin{table}
$\begin{array}{|r|c|cccc|c|r|c|cccc|}
\cline{1-6} \cline{8-13}
\text{gp} & K & a^{K-3}_{\mathfrak{m}} & a^{K-2}_{\mathfrak{m}}
& a^{K-1}_{\mathfrak{m}} & a^K_{\mathfrak{m}} & \vphantom{\Bigl(}\quad &
\text{gp} & K & a^{K-3}_{\mathfrak{m}} & a^{K-2}_{\mathfrak{m}}
& a^{K-1}_{\mathfrak{m}} & a^K_{\mathfrak{m}} \\
\cline{1-6} \cline{8-13}
36 & 5 & -5 & -5 & -5 & -5 & &
850 & 5 & -\infty & -6 & -5 & -5 \\
48 & 5 & -5 & -5 & -5 & -5 & &
852 & 4 & -\infty & -5 & -4 & -4 \\
52 & 4 & -4 & -4 & -4 & -4 & &
853 & 4 & -\infty & -5 & -4 & -4 \\
194 & 5 & -6 & -6 & -5 & -5 & &
854 & 4 & -\infty & -5 & -4 & -4 \\
513 & 5 & -\infty & -5 & -5 & -5 & &
859 & 4 & -\infty & -6 & -4 & -4 \\
515 & 5 & -5 & -5 & -5 & -5 & &
860 & 4 & -\infty & -5 & -4 & -4 \\
527 & 4 & -\infty & -4 & -4 & -4 & &
866 & 4 & -\infty & -5 & -4 & -4 \\
551 & 5 & -6 & -5 & -5 & -5 & &
928 & 4 & -\infty & -\infty & -4 & -4 \\
560 & 4 & -6 & -5 & -4 & -4 & &
929 & 4 & -\infty & -5 & -4 & -4 \\
561 & 4 & -5 & -4 & -4 & -4 & &
931 & 4 & -\infty & -3 & -4 & -4 \\
621 & 5 & -\infty & -5 & -5 & -5 & &
932 & 4 & -\infty & -6 & -4 & -4 \\
623 & 4 & -\infty & -4 & -4 & -4 & &
934 & 4 & -\infty & -3 & -5 & -4 \\
630 & 5 & -\infty & -5 & -5 & -5 & &
1578 & 6 & -\infty & -6 & -6 & -6 \\
635 & 4 & -\infty & -4 & -4 & -4 & &
1615 & 4 & -\infty & -\infty & -4 & -4 \\
636 & 4 & -\infty & -4 & -4 & -4 & &
1620 & 4 & -\infty & -4 & -4 & -4 \\
642 & 4 & -\infty & -4 & -4 & -4 & &
1735 & 5 & -\infty & -5 & -5 & -5 \\
643 & 4 & -\infty & -4 & -4 & -4 & &
1751 & 4 & -\infty & -6 & -4 & -4 \\
645 & 4 & -\infty & -\infty & -4 & -4 & &
1753 & 4 & -\infty & -\infty & -4 & -4 \\
646 & 4 & -\infty & -4 & -4 & -4 & &
1755 & 5 & -\infty & -\infty & -5 & -5 \\
740 & 4 & -\infty & -4 & -4 & -4 & &
1757 & 4 & -\infty & -\infty & -4 & -4 \\
742 & 4 & -\infty & -4 & -4 & -4 & &
1758 & 4 & -\infty & -\infty & -5 & -4 \\
753 & 5 & -\infty & -5 & -5 & -5 & &
1759 & 4 & -\infty & -\infty & -4 & -4 \\
761 & 5 & -5 & -5 & -5 & -5 & &
1800 & 4 & -\infty & -4 & -4 & -4 \\
764 & 4 & -\infty & -5 & -4 & -4 & &
2216 & 5 & -\infty & -\infty & -5 & -5 \\
780 & 4 & -6 & -4 & -4 & -4 & &
2222 & 5 & -\infty & -6 & -5 & -5 \\
801 & 4 & -4 & -4 & -4 & -4 & &
2264 & 5 & -\infty & -6 & -5 & -5 \\
813 & 4 & -4 & -4 & -4 & -4 & &
2317 & 4 & -\infty & -\infty & -4 & -4 \\
823 & 4 & -4 & -4 & -4 & -4 & &
2326 & 4 & -\infty & -\infty & -\infty & -4 \\
\cline{8-13}
836 & 4 & -4 & -5 & -4 & -4 \\
\cline{1-6}
\end{array}$
\vspace*{5pt}
\caption{For each of the 57 groups of order $128$ with $\gtD(G)=3$, we give
its number in the Small Groups library, the Krull dimension $K$
and the last four $a$-invariants of $H^*(G)$. For
$0 \leq i < K-3$ one has $a^i_{\mathfrak{m}} = -\infty$.
The defect $\CMd(G)$ is three if and only if
$a^{K-3}_{\mathfrak{m}}$ is finite.}
\label{table:128aInvts}
\end{table}

In order to calculate the $a$-invariants, one uses Lemma~4.3
of~\cite{Benson:DicksonCompCoho}, which is based on methods of N.~V. Trung\@.
The lemma says that if $\zeta \in R^n$ is such that $\Ann_M(\zeta)$ consists
entirely of $\mathfrak{m}$-torsion, then
\begin{equation}
\label{eqn:Benson}
a^{i+1}_{\mathfrak{m}}(M) + n \leq a^i_{\mathfrak{m}}(M/\zeta M)
\leq \max(a^i_{\mathfrak{m}}(M), a^{i+1}_{\mathfrak{m}}(M) + n) \, .
\end{equation}
Replacing $\zeta$ by a suitable power if necessary, one can
arrange for $a^{i+1}_{\mathfrak{m}}(M) + n \geq a^i_{\mathfrak{m}}(M)$
and therefore $a^{i+1}_{\mathfrak{m}}(M) = a^i_{\mathfrak{m}}(M/\zeta M) - n$.
So the $a$-invariants of $H^*(G)$ may be computed recursively.
To start the recursion we need
$a^0_{\mathfrak{m}}(M)$, which is the $\mathfrak{m}$-torsion:
so if $\zeta \in R^n$ is such that $\Ann_M(\zeta)$ is finite dimensional,
then $a^0_{\mathfrak{m}}(M)$ is the top dimension of $\Ann_M(\zeta)$.

Hence by starting from a filter-regular system of parameters and
raising some of the parameters to higher powers if necessary, one
may compute the $a$-invariants of $H^*(G)$ by just computing kernels.
For a cohomology computation one chooses parameters in low degrees. So it is
perhaps surprising that a survey of the author's computations of all
256 nonabelian groups of order 64 and some 61 groups of order 128 led to
precisely one case where powers of the chosen parameters were necessary.
This is the Sylow $2$-subgroup of $L_3(4)$ which has $a$-invariants
$-\infty,-\infty,-3,-5,-4$: a filter-regular system of parameters
in degrees $4,4,2,2$ led to kernels with top degrees $-\infty,-\infty,5,5,8$,
leading to problems with the calculation of $a^3_{\mathfrak{m}}$. Squaring
the third parameter led to kernels with top degrees
$-\infty,-\infty,5,7,10$, which was sufficient to permit calculation of
the $a$-invariants.

\section{Excess and defect}

\begin{defn}
As in the introduction let $G$ be a finite group, $p$ a prime number and
$k$ a field of characteristic~$p$. We define
the Duflot excess $\De(G)=\De_p(G)$ by
\[
\De_p(G) = \depth H^*(G,k) - \prank(Z(S)) \, .
\]
\end{defn}

\noindent
The following inequalities follow immediately from this definition taken
together with Equations
\eqref{eqn:QuillenDuflot}~and \eqref{eqn:gtD-CMd}\@.
\begin{xalignat}{3}
\label{eqn:excess}
0 \leq \De(G) & \leq \gtD(G) &
\CMd(G) + \De(G) & = \gtD(G) &
\De(G) & \geq \De(S) \, .
\end{xalignat}

Quillen showed that the extraspecial $2$-group $G = 2^{1+2n}_+$ has
Cohen--Macaulay cohomology~\cite{Quillen:Extraspecial}.
So this group has $\De(G)=n$ and $\CMd(G)=0$.

Now let $p$ be an odd prime, and let $G$ be the extraspecial $p$-group
$G=p^{1+2n}_+$ of exponent~$p$.  With the single exception of the case
$(p,n)=(3,1)$, Minh proved~\cite{Minh:EssExtra} that this group has
$\CMd(G)=n$ and $\De(G)=0$. In the one exceptional case the cohomology ring
is Cohen--Macaulay~\cite{MilgramTezuka}\@.

One good way to produce groups with small $\De(G)/\CMd(G)$ ratio
satisfying Conjecture~\ref{conj:Reg}
is by iterating the wreath product construction. By passing
from $H$~to $H \wr C_p$ one multiplies the $p$-rank by~$p$ but increases the
depth by one only~\cite{CaHe:Wreath}\@.

\begin{ques}
How (for large values of~$n$) are the $p$-groups of order $p^n$ distributed
on the graph with $\CMd(G)$ on the $x$-axis and $\De(G)$ on the $y$-axis?
\end{ques}

\section{Outlook}
\noindent
To test the conjecture further we need to find more high defect groups.

There are 24 groups of order~$3^6$ with $\gtD(G)=3$. The presence of essential
classes in low degrees demonstrates that at least three of these groups
have $\CMd(G)=3$. These groups are numbers $35$, $56$ and $67$ in the Small
Groups Library. There are essential classes in degrees $4$, $2$ and $4$
respectively.
Recall from Carlson's paper~\cite{Carlson:DepthTransfer} that the presence
of essential classes means that $\depth H^*(G,k) = \prank Z(G)$ and
therefore $\CMd(G)=\gtD(G)$.

Group number $299$ of order $256$ has $\gtD(G)=4$. The presence of an essential
class in $H^3(G)$ means that $\CMd(G)=4$ too.

\appendix
\section{Computing the $p$-rank}
\label{sect:appendix}
The Small Groups library was accessed from GAP~\cite{GAP4}\@. There is
no built-in command in GAP that returns the $p$-rank of a given $p$-group.
The simplest way to calculate it using existing commands would to be to
generate the entire subgroup lattice and then filter out the elementary
abelian subgroups. We chose instead to enumerate the conjugacy classes
of elementary abelian subgroups using a straightforward if not particularly
efficient inductive approach.

Let $G$ be a $p$-group. If $\left|G\right|$ is small then it is feasible to
list all the elements of the group. By testing each element one then obtains
the list of all order~$p$ elements. A further element by element test yields
all central elements of order~$p$. This is one way
to obtain the greatest central elementary abelian subgroup $\Omega_1(Z(G))$,
denoted~$C$ in the paper.
Of course, if one only wants $\Omega_1(Z(G))$ then it is quicker to make
use of the
function \texttt{IndependentGeneratorsOfAbelianGroup}\@.

Carlson~\cite{Carlson:DepthTransfer} shows that one only needs
the elementary abelians which contain $\Omega_1(Z(G))$. Given such an
elementary abelian $V$ of order~$p^d$, one can list all the order $p$ elements
in $C_G(V)$ and so obtain all the order $p^{d+1}$ elementary abelians
containing~$V$. This is the inductive step: the induction starts with
$V=C$.

\bibliographystyle{abbrv}
\bibliography{../united}



\end{document}